\documentclass[11pt]{amsart}   	
\usepackage{geometry}                		
\usepackage{graphicx}				
\usepackage{amssymb}
\usepackage{amsmath}
\usepackage{enumerate}

\newcommand{\cS}{{\mathcal S}}
\newcommand{\cN}{{\mathcal N}}

\newcommand{\cT}{{\mathcal T}}
\newcommand{\cU}{{\mathcal U}}

\usepackage{color}
\newcommand{\blue}{\textcolor{black}}
\newtheorem{theorem}{Theorem}[section]

\title[Universal tree-based network]{A universal tree-based network with the minimum number of reticulations}

\author{Magnus Bordewich}
\address{School of Engineering Computer Sciences, Durham University, Durham DH1 3LE, United Kingdom}
\email{m.j.r.bordewich@durham.ac.uk}

\author{Charles Semple}
\address{School of Mathematics and Statistics, University of Canterbury, Christchurch, New Zealand}
\email{charles.semple@canterbury.ac.nz}



\date{\today}

\begin{document}

\maketitle

\begin{abstract}
A tree-based network $\cN$ on $X$ is universal if every rooted binary phylogenetic $X$-tree is a base tree for $\cN$. \blue{Hayamizu} and, independently, Zhang constructively showed that, for all positive integers $n$, there exists an universal tree-based network on $n$ leaves. For all $n$, \blue{Hayamizu's} construction contains $\Theta(n!)$ reticulations, while Zhang's construction contains $\Theta(n^2)$ reticulations. A simple counting argument shows that an universal tree-based network has $\Omega(n\log n)$ reticulations. With this in mind, \blue{Hayamizu as well as} Steel posed the problem of determining whether or not such networks exists with $O(n\log n)$ reticulations. In this paper, we show that, for all $n$, there exists an universal tree-based network on $n$ leaves with $O(n\log n)$ reticulations.\end{abstract}

\section{Introduction}

A {\em phylogenetic network $\cN$ on $X$} is a rooted acyclic digraph with no edges in parallel satisfying the following properties:
\begin{enumerate}[(i)]
\item the root has out-degree two;
\item a vertex with out-degree zero has in-degree one, and the set of vertices with out-degree zero is $X$; and
\item all other vertices either have in-degree one and out-degree two, or in-degree two and out-degree one.
\end{enumerate}
For technical reasons, if $|X|=1$, then we additionally allow $\cN$ to consist of the single vertex in $X$. As described here, a phylogenetic network is sometimes referred to as a binary phylogenetic network. Vertices of in-degree one and out-degree zero are called {\em leaves}, while vertices of in-degree one and out-degree two are {\em tree vertices} and vertices of in-degree two and out-degree one are {\em reticulations}. \blue{An edge directed into a reticulation is called a {\em reticulation edge}. All other edges are {\em tree edges}.} A {\em rooted binary phylogenetic $X$-tree} is a phylogenetic network on $X$ with no reticulations. Throughout the paper, we denote the size of $X$ by $n$.

Let $\cT$ be a phylogenetic $X$-tree and let $\cN$ be a phylogenetic network on $X$. We say that $\cN$ {\em displays} $\cT$ if, up to suppressing degree-two vertices, $\cT$ can be obtained from $\cN$ by deleting edges and vertices, in which case, the resulting acyclic digraph together with a path to its root from the root of $\cN$ is an {\em embedding} of $\cT$ in $\cN$. Note that if $\cS$ is an embedding of $\cT$ in $\cN$, then the root of $\cS$ is the root of $\cN$ and so it may have out-degree one.

A phylogenetic network $\cN$ on $X$ is a {\em tree-based network} if there is an embedding $\cS$ of a phylogenetic $X$-tree $\cT$ in $\cN$ such that $\cS$ contains every vertex of $\cN$. If this holds, then $\cS$, as well as $\cT$, is a {\em base tree} for $\cN$. Tree-based networks were introduced by Francis and Steel~\cite{fra15} as a way of quantifying the concept of an `underlying tree'. In particular, tree-based networks can be equivalently defined as those phylogenetic networks that can be obtained from a rooted binary phylogenetic tree $\cT$ by simply adding edges whose end-vertices subdivided edges of $\cT$. The concept of a tree-based network is relevant to the on-going debate concerning the extent to which the evolution of early life can be viewed as simply a phylogenetic tree with some additional edges or whether the concept of underlying tree is completely meaningless~\cite{abb12, dag06}. Tree-based networks have generated considerable interest in the last year (for example, see~\cite{ana16, hay16, jet17, sem16, zha16}).

A tree-based network on $X$ is {\em universal} if every rooted binary phylogenetic $X$-tree is a base tree for $\cN$. Not all phylogenetic networks are tree-based and so, \textit{a priori}, it is not clear whether a universal tree-based network on $X$ exists for all positive integers $n$. Francis and Steel~\cite{fra15} showed that if $|X|\le 3$, then there exists a universal tree-based network on $X$, and posed the problem of determining whether or not such a network exists for all $n$. \blue{Hayamizu~\cite{hay16}} and, independently, Zhang~\cite{zha16} constructively showed that, for all $n$, there is indeed a universal tree-based network on $n$ leaves. For all $n$, the construction in~\cite{hay16} contains $\Theta(n!)$ reticulations, while the construction in~\cite{zha16} contains $\Theta(n^2)$ reticulations. As a consequence of this, \blue{Hayamizu~\cite{hay16} and} Mike Steel (private communication) asked the very natural question: What is the minimum number of reticulations in a universal tree-based network?

Let $b_n$ denote the number of rooted binary phylogenetic $X$-trees. A classical result in phylogenetics dates back to Schr\"{o}der (1870) who showed that
$$b_n=1\times 3\times 5\times \cdots\times (2n-3)=\frac{(2n-2)!}{(n-1)!2^{n-1}}.$$
Therefore, by Stirling's approximation,
\begin{align}
\label{stirling}
b_n\sim \frac{1}{\sqrt{2}}\left(\frac{2}{e}\right)^n n^{n-1}.
\end{align}
Now, if $\cN$ is a phylogenetic network on $X$ with $r$ reticulations, then $\cN$ displays at most $2^r$ distinct rooted binary phylogenetic trees, as each embedding of such a tree is realised by choosing exactly one of the reticulation edges directed into each reticulation. Equating this with (\ref{stirling}) and solving for $r$, it follows that a phylogenetic network on $X$ that displays every rooted binary phylogenetic $X$-tree has $\Omega(n\log n)$ reticulations. \blue{Hayamizu's and} Steel's question is therefore more precisely stated as the following: for all positive integers $n$, does there exist an universal tree-based network on $n$ leaves with $O(n\log n)$ reticulations? In this paper, we affirmatively answer this question. In particular, we establish the following theorem, the proof of which is given in the next section. Note that all logarithms in this paper are base~$2$.

\begin{theorem}
For all positive integers $n$, there exists a universal tree-based network on $n$ leaves with $O(n\log n)$ reticulations.
\label{main}
\end{theorem}

We end the introduction with \blue{two remarks}. The constructions in~\cite{hay16} and~\cite{zha16} are very similar. Intuitively, for each $n$, they construct a tree-based network consisting of two halves. The first half, which contains the root, embeds all tree shapes, that is, all (unlabelled) rooted binary trees. The second half, which contains the leaves, then reorders the leaves of these trees to produce the desired rooted binary phylogenetic tree. The first half constructions of both papers are identical and consists of $\Theta(n^2)$ reticulations. The difference between the two constructions lies in the second halves. Although the approaches are the same, Zhang~\cite{zha16} uses $\Theta(n^2)$ reticulations, while \blue{Hayamizu~\cite{hay16}} use $\Theta(n!)$ reticulations. We adopt a similar overall approach in this paper, except that both our first and second half constructions use $O(n\log n)$ reticulations.

\blue{The second remark is that, like the universal tree-based networks in~\cite{hay16} and~\cite{zha16}, the universal tree-based network constructed in this paper is both temporal and stack-free. A phylogenetic network $\cN$ is {\em temporal} if there is a mapping from the set of vertices of $\cN$ to the non-negative integers such that if $(u, v)$ is a tree edge, then $t(u) < t(v)$, while if $(u, v)$ is a reticulation edge, then $t(u)=t(v)$. Biologically, if $\cN$ is temporal, then $\cN$ satisfies two natural timing constraints. In particular, successively occurring speciation events, and contemporaneously occurring reticulation events. Lastly, a phylogenetic network is {\em stack-free} if it has no two reticulations one of which is a parent of the other.}

\section{Proof of Theorem~\ref{main}}

\begin{figure}
\begin{center}
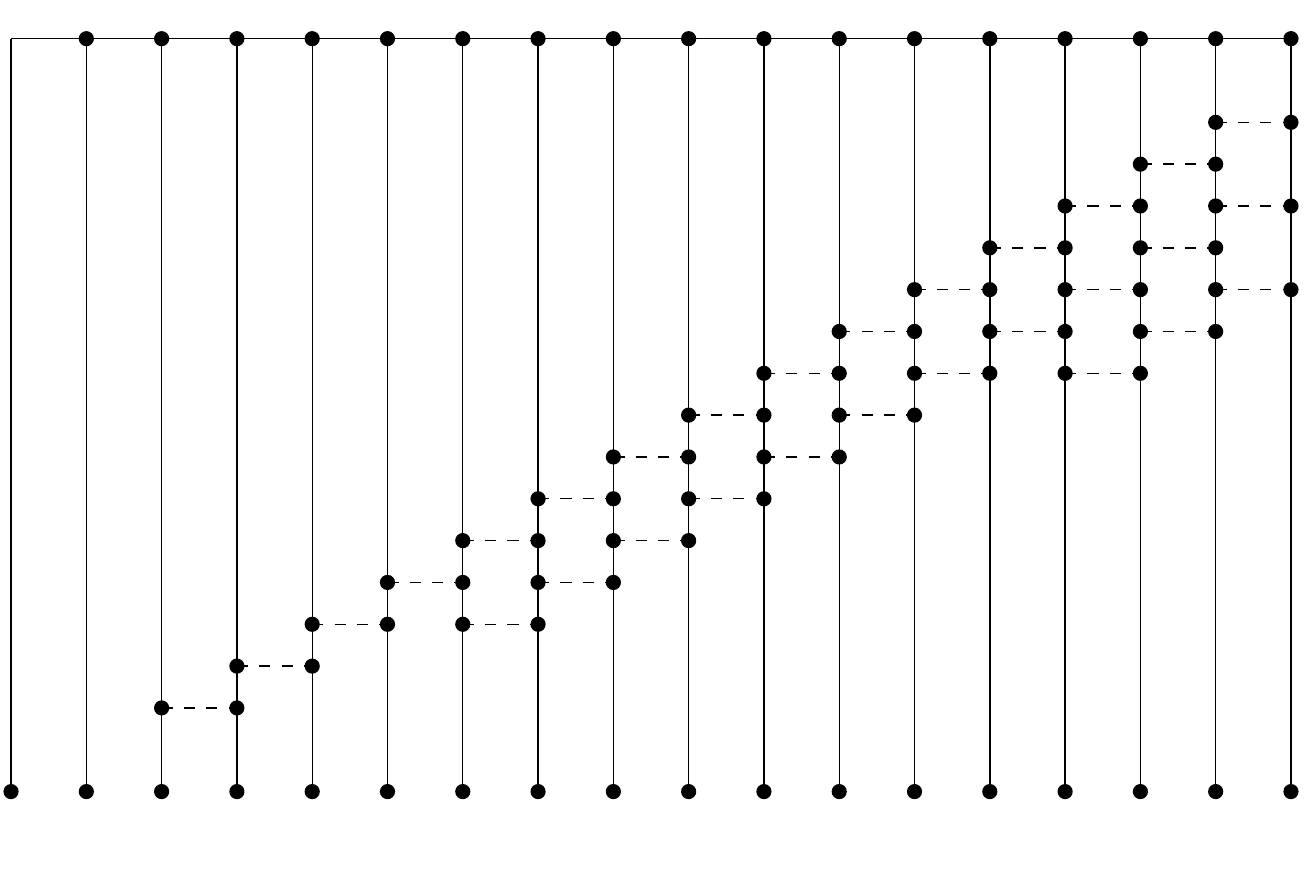
\caption{The construction $A_{18}$, where the root is denoted by $\rho$. The rooted binary caterpillar $(1, 2, \ldots, 18)$ is shown as solid edges, while dashed edges are the additional edges. Vertical edges are directed down the page, and horizontal edges are directed right to left.}
\label{Un}
\end{center}
\end{figure}

\begin{figure}
\begin{center}
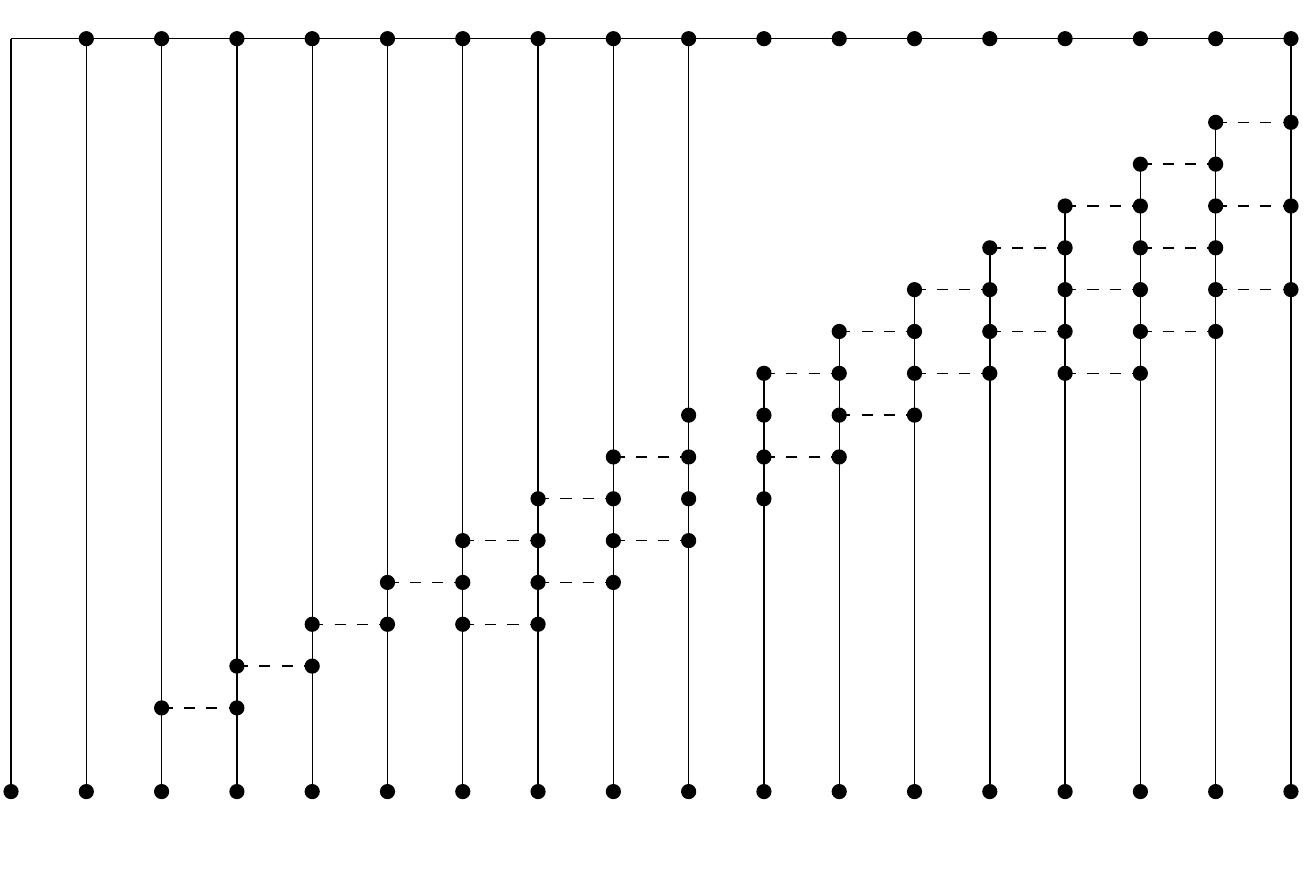
\caption{The network $A_{18, 10}$, where the root is denoted by $\rho$. Vertical edges are directed down the page, and horizontal edges are directed from right to left.}
\label{split}
\end{center}
\end{figure}

For all positive integers $n$, we first construct a phylogenetic network $\cU_n$, which we eventually establish is an universal tree-based network on $n$ leaves with $O(n\log n)$ reticulations. We begin by describing the structure of the top half of $\cU_n$, which we denote by $A_n$. The {\em rooted binary caterpillar $(1, 2, \ldots, n)$} is the rooted binary phylogenetic tree on $\{1, 2, \ldots, n\}$ such that the parents of leaves $1$ and $2$ are the same, the parent of leaf $n$ is the root, and $q_n, q_{n-1}, \ldots, q_2, 1$ is a directed path from the root to leaf $1$ where, for each $i\in \{2, 3, \ldots, n\}$, $q_i$ is the parent of leaf $i$ (see Fig.~\ref{Un}, solid edges). The path $q_n, q_{n-1}, \ldots, q_2, 1$ is the {\em spine} of the caterpillar. The digraph $A_n$ can be viewed as the rooted binary caterpillar $(1, 2, \ldots, n)$ with additional edges whose end-vertices subdivide `neighbouring' pendant edges of the caterpillar. More precisely, take the rooted binary caterpillar $(1, 2, \ldots, n)$ and, for each $i\geq 2$, subdivide the pendant edge incident with leaf $i$ with
$$\left\lfloor\log \left(\frac{i}{2}\right)\right\rfloor + \left\lfloor\log \left(\frac{i+1}{2}\right)\right\rfloor$$ 
new vertices. New edges are now added between the new vertices as follows. For each $i$, denote the path from the spine to leaf $i$ by $P_i$. Numbering the new vertices on each $P_i$ starting at the spine, add an edge directed from the $2j$-th new vertex on $P_i$ to the $(2j-1)$-th new vertex on $P_{i-1}$. \blue{Observe that the number of edges joining $P_i$ and $P_{i-1}$ and the placement of these edges is determined by $i$ and not $n$.} For an example, see Fig.~\ref{Un} which shows the resulting construction for when $n=18$ where, as also in Fig.~\ref{split}, vertical edges are directed down the page, and horizontal edges are directed from right to left. Note that, for each $i$, the odd numbered new vertices on $P_i$ are reticulations, and the even numbered new vertices on $P_i$ are tree vertices. The total number of reticulations in $A_n$ is
$$\sum_{i=2}^n \left\lfloor\log \left(\frac{i+1}{2}\right)\right\rfloor\leq \log \left(\prod_{i=2}^n \left(\frac{i+1}{2}\right)\right)\leq \log \left(n^n\right) = n\log n.$$
Ignoring the leaf labels, we will later show that if $T$ is a rooted binary tree, then there is an embedding of $T$ in $A_n$ using every vertex.

\begin{figure}
\center
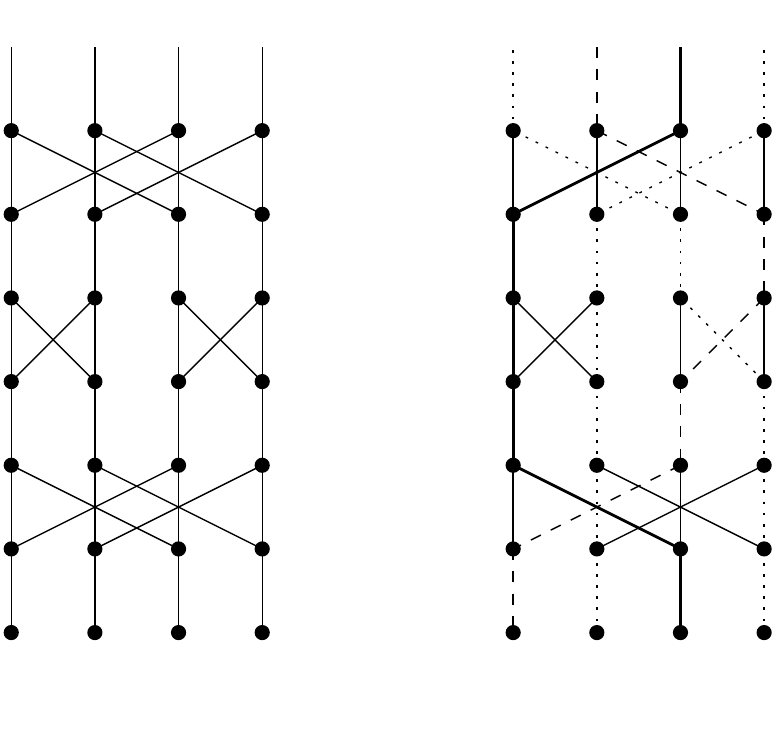
\caption{\blue{(i) The construction $B_4$. The vertical edges and the six pairs of ``crossover'' edges correspond to the four wires and the six switches of the Bene\v{s} network of size four, respectively. (ii) An embedded reordering of the permutation $(4, 1, 3, 2)$. In (i) and (ii), edges are directed down the page.}}
\label{benes}
\end{figure}
\begin{figure}
\center
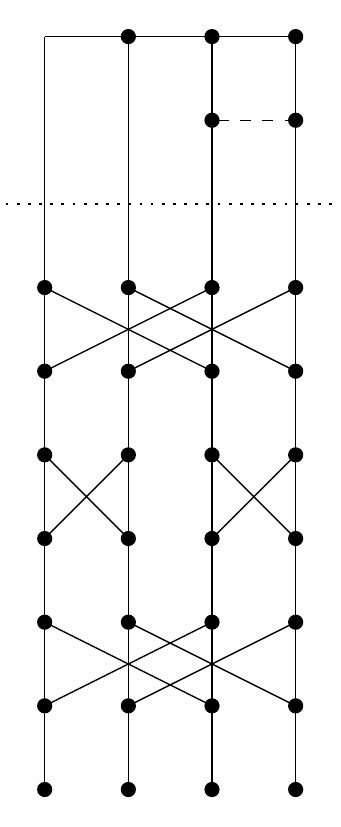
\caption{\blue{Illustrating how $A_4$ and $B_4$ combine to give $\cU_4$. The dotted line separates the two halves. Horizontal edges are directed right to left, while all other edges are directed down the page.}}
\label{twohalves}
\end{figure}

\blue{The bottom half of $\cU_n$, denoted $B_n$, uses a {\em Bene\v{s} network} construction to enable any permutation of the leaves. A Bene\v{s} network~\cite{ben64a, ben64b} is an example of a ``rearrangeable non-blocking'' network. In particular, it is an arrangement of wires (transmission links) and switches which realises all possible permutations between the input and output terminals. Each switch corresponds to a permutation of size two. The original motivation for such networks is in providing a telephone service in which it is possible to rearrange existing calls to allow any new call into the system. For further background on Bene\v{s} networks, we refer the interested reader to~\cite{lei92}. For the construction of $B_n$, we extend, for each $i$, the path $P_i$ of $A_n$, each extension corresponds to a wire and, whenever there is a switch between wires $i$ and $j$ in the Bene\v{s} network, we insert a tree vertex followed by a reticulation in each corresponding path, making the new tree vertices on paths $P_i$ and $P_j$ the parents of the new reticulations on paths $P_j$ and $P_i$, respectively. To illustrate, Fig.~\ref{benes} shows (i) the construction of $B_4$ from the Bene\v{s} network on four wires, and (ii) the embedding of a permutation that reorders the leaves $(4, 1, 3, 2)$ coming out of $A_4$ to the ordering $(1, 2, 3, 4)$. Note that, for each ordering of the elements in $\{1, 2, 3, 4\}$, there is an embedding of the permutation that reorders it to the ordering $(1, 2, 3, 4)$. Furthermore, for $n=4$, Fig.~\ref{twohalves} shows how the two halves are combined.}


\begin{figure}
\center
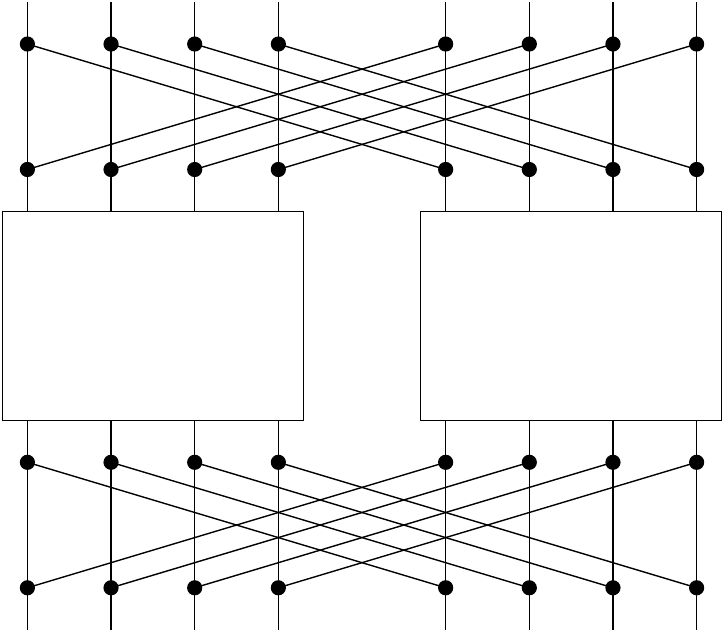
\caption{\blue{The Bene\v{s} network of size~$8$ constructed from two Bene\v{s} networks of size~$4$. Edges are directed down the page.}}
\label{benes8}
\end{figure}

\blue{Bene\v{s} networks were originally constructed for when $n$ is a power of two. However, more recently, Bene\v{s} networks have been constructed for arbitrary $n$~\cite{cha97}. Both constructions are recursive in that the Bene\v{s} network of size $n$ (that is, with $n$ input and $n$ output terminals) is built from the Bene\v{s} network of size $\left\lfloor\frac{n}{2}\right\rfloor$ and the Bene\v{s} network of size $\left\lceil\frac{n}{2}\right\rceil$. For example, the Bene\v{s} network of size~$8$ is constructed from two Bene\v{s} networks of size~$4$ as shown in Fig.~\ref{benes8}, where each switch has been replaced with a pair of crossover edges. Moreover, both when $n$ is a power of two and in the general case, the constructions have $O(n\log n)$ switches. Thus, our network $B_n$ has $O(n\log n)$ reticulations and so, by construction, for all $n$, the phylogenetic network $\cU_n$ has $O(n\log n)$ reticulations. Furthermore, for all $n$, we have that $\cU_n$ is tree-based since we can take the vertical edges (that is, the caterpillar in $A_n$ and the identity permutation in $B_n$) as a base tree. Since the embedding of any permutation contains every vertex of $B_n$, to complete the proof of Theorem~\ref{main} it suffices to show that, up to leaf labels, every rooted binary phylogenetic tree on $\{1, 2, \ldots, n\}$ is a base tree of $A_n$.}

We use induction on $n$, to show that $A_n$ embeds any tree shape on $n$ leaves and that every vertex of $A_n$ is contained in the embedding. Since there is exactly one tree shape on $n\leq 3$ leaves, the base case holds for all $n\leq 3$. Now suppose that $n\ge 4$ and that, for all $m\le n-1$, the construction $A_m$ embeds any tree shape with $m$ leaves and every vertex of $A_m$ is contained in the embedding. 

Let $T$ be a tree shape on $n$ leaves, and denote the two maximal rooted binary subtrees whose roots are children of the root of \blue{$T$} by  $T_1$ and $T_2$. Let $t_1=|T_1|$ and $t_2=|T_2|$. Without loss of generality we may assume that $t_1\geq t_2$. Note that $t_2\leq n/2$. Consider the network, denoted $A_{n, t_1}$, obtained from $A_n$ by removing all of the additional edges added between $P_{t_1}$ and $P_{t_1+1}$ in the construction of $A_n$, and also removing the first edge on each of the paths $P_{t_1+1}, P_{t_1+2}, \ldots, P_{n-1}$. The network $A_{18, 10}$ is shown in Fig.~\ref{split}.

Up to degree-two vertices, the root of $A_{n, t_1}$ has two children. \blue{The first is the root of a subnetwork $D_{t_1}$ which is isomorphic to $A_{t_1}$ since, in the construction, the number of edges between $P_i$ and $P_{i-1}$ depend on $i$ only and not $n$}. The second is the root of a subnetwork $D_{t_2}$. Now, up to isomorphism, $A_{t_2}$ can be obtained from $D_{t_2}$ by deleting a subset of the additional edges which were added to the rooted binary caterpillar $(1, 2, \ldots, n)$ to create $A_{n}$, and suppressing degree-two vertices. To see this and, in particular, that we have enough `additional' edges, observe that, as $t_2\le n/2$, the $l$-th leaf of $D_{t_2}$ corresponds to at least the $2l$-th leaf of $A_n$. Hence, there are at least $$\left\lfloor\log \left(\frac{2(l+1)}{2}\right)\right\rfloor-1=\left\lfloor\log\left(\frac{l+1}{2}\right)\right\rfloor$$ edges between the paths $P'_l$ and $P'_{l+1}$ of $D_{t_2}$, where, for all $l\in \{t_1+1, t_2+2, \ldots, n-1\}$, the path $P'_i$ is the path obtained from $P_i$ (in $A_n$) by deleting the first edge. Thus we can delete edges of $D_{t_2}$ starting with those nearest the leaves until we obtain the correct number of edges between $P'_l$ and $P'_{l+1}$ for $A_{t_2}$. Finally, we can suppress the \blue{resulting} degree-two vertices which \blue{were created in one of two ways. Firstly,} by the initial deletion of edges adjacent to the spine, and are therefore now on the unique path from the root to leaf $t_1+1$ \blue{and, secondly,} created in the later deletions, and are therefore on paths of degree-two vertices between a leaf and \blue{its} first ancestor which has degree greater than two (thus on every path from the root to that leaf). \blue{Hence,} if we take any embedding of $T_2$ in $A_{t_2}$ which contains every vertex of $A_{t_2}$, it corresponds to an embedding of $T_2$ in $D_{t_2}$ \blue{which} contains every vertex of $D_{t_2}$.

By induction, $D_{t_1}$ and $D_{t_2}$ contain an embedding of $T_1$ and $T_2$ that uses each of its vertices, respectively. Furthermore, extending these embeddings in $A_n$ by including the vertices in the unique paths from the root of $A_n$ to the roots of $D_{t_1}$ and $D_{t_2}$ gives an embedding of $T$ in $A_n$ that contains every vertex of $A_n$. This completes the proof of Theorem~\ref{main}.

\section*{Acknowledgements} \blue{We thank one of the anonymous referees for suggesting the simpler Bene\v{s} networks as an alternative to sorting networks for constructing the bottom half of $\cU_n$.}

\end{document}